\documentclass[lettersize,journal]{IEEEtran}

\usepackage{amsmath,amsfonts,amssymb}
\usepackage{mathtools}
\usepackage{bm}
\usepackage{breqn}
\usepackage{hyperref}
\usepackage{siunitx}
\DeclareSIUnit{\EUR}{\text{\euro}}

\usepackage[dvipsnames]{xcolor}
\usepackage{cleveref}

\hypersetup{
    colorlinks=true,
    linkcolor=blue,
    citecolor=ForestGreen,
    urlcolor=blue
}

\crefformat{algorithm}{#2{\color{BrickRed}#1}#3}
\Crefformat{algorithm}{#2{\color{BrickRed}#1}#3}
\usepackage{mathrsfs}

\usepackage{amsthm}
\theoremstyle{plain}%

\usepackage{rotating}
\usepackage{longtable}
\usepackage{booktabs}

\usepackage{eurosym}
\usepackage{subcaption}

\usepackage{pgfplots}
\pgfplotsset{compat=1.18}
\usepackage{xcolor}

\usepackage{pgfplotstable}
\usepackage{xstring}

\usepackage{tikz}
\usetikzlibrary{patterns,calc,arrows.meta,positioning}
\makeatletter
\tikzset{
  hatch distance/.store in=\hatchdistance,
  hatch distance=5pt,
  hatch thickness/.store in=\hatchthickness,
  hatch thickness=1.2pt
}
\pgfdeclarepatternformonly[\hatchdistance,\hatchthickness]{custom hatch}
  {\pgfqpoint{-1pt}{-1pt}}{\pgfqpoint{\hatchdistance}{\hatchdistance}}
  {\pgfpoint{\hatchdistance}{\hatchdistance}}{
    \pgfsetlinewidth{\hatchthickness}
    \pgfpathmoveto{\pgfqpoint{0pt}{0pt}}
    \pgfpathlineto{\pgfqpoint{\hatchdistance}{\hatchdistance}}
    \pgfusepath{stroke}
}
\makeatother

\makeatletter
\let\@msm@th@eqref\eqref
\renewcommand{\eqref}[1]{%
  \begingroup
  \leavevmode
  \color{blue}%
  \hypersetup{linkbordercolor=[named]{blue}}%
  \@msm@th@eqref{#1}%
  \endgroup
}
\makeatother

\usepackage{setspace}

\begin{document}

\title{Towards time-variant scenario reduction for energy system optimization modeling under uncertainty}

% \author{Y. Werner, J. M. Morales,~\IEEEmembership{Senior Member,~IEEE}, S. Pineda,~\IEEEmembership{Senior Member,~IEEE}, and S. Wogrin,~\IEEEmembership{Senior Member,~IEEE} \vspace{-0.4cm} %
\author{Y. Werner, J. M. Morales, S. Pineda, and S. Wogrin \vspace{-0.4cm} %
\thanks{
Y. Werner and S. Wogrin are with the Institute of Electricity Economics and Energy Innovation at Graz University of Technology, 8020 Graz, Austria.\\
E-mails: $\{$yannick.werner, wogrin$\}$@tugraz.at.\\
J. M. Morales and S. Pineda are with the research group OASYS, University of Málaga, Málaga 29071, Spain. E-mails: $\{$spineda,juan.morales$\}$@uma.es.\
}
}



\vspace{-1cm}
\IEEEaftertitletext{\vspace{-0.6\baselineskip}}
\maketitle

% \newcommand\copyrighttext{%
%   \footnotesize \textcopyright \the\year{} IEEE. Personal use of this material is permitted. Permission from IEEE must be obtained for all other uses, including reprinting/republishing this material for advertising or promotional purposes, collecting new collected works for resale or redistribution to servers or lists, or reuse of any copyrighted component of this work in other works.}

% \newcommand\copyrightnotice{%
% \begin{tikzpicture}[remember picture,overlay]
% \node[anchor=south,yshift=10pt] at (current page.south) {\fbox{\parbox{\dimexpr0.75\textwidth-\fboxsep-\fboxrule\relax}{\copyrighttext}}};
% \end{tikzpicture}%
% }

% \newcommand\IEEEcopyrighttext{\textbf{979-8-3195-3554-2/26/\$31.00 ©2026 IEEE}}
% \newcommand\IEEEcopyrightnotice{%
% \begin{tikzpicture}[remember picture,overlay]
% \node[anchor=south west, xshift=1.6cm, yshift=0.0cm] 
% at (current page.south west) {%
%  \IEEEcopyrighttext};
% \end{tikzpicture}%
% }

%for ARXIV use this command
% \copyrightnotice 
% for the "real" upload you can use that:
%\IEEEcopyrightnotice
% \vspace{3pt}

\begin{abstract}
Stochastic programming has become a popular tool for supporting decision-making under uncertainty in the long-term planning of energy systems. Existing scenario reduction methods, however, are naive about the long-term temporal nature of scenarios, which limits their efficiency in reducing model size. In this paper, we overcome this inefficiency by proposing a novel time-variant scenario reduction framework that explicitly allows for varying scenario aggregations over time. As a result, scenario probabilities become time-variant, enabling not only the accurate capture of scenario realizations but also their probabilities at the time steps that drive investment decisions. This substantially increases flexibility compared to traditional time-invariant methods, which we demonstrate on a two-stage stochastic generation expansion planning problem with uncertain renewable power production.
\end{abstract}


\section{Introduction}
\label{sec:intro}

Long-term two-stage stochastic energy system optimization models (ESOMs), such as AnyMOD \cite{Goeke2024}, EMPIRE \cite{Backe2022}, and PyPSA~\cite{Brown2018a}, have become important tools for supporting decision-making under uncertainty in energy system planning. They typically optimize investments in the first stage and system operation in the second stage, accounting for uncertainty, e.g., from intermittent renewable power production. When spanning time horizons of 20–60 years and considering 10--40 scenarios to represent the uncertainty, the model size grows rapidly. 
Scenario reduction techniques aim to overcome this by reducing the number of scenarios required to accurately represent the original distribution of the random variables without distorting the optimal solution of the stochastic optimization problem. Pioneering works, such as those by Dupačová et al.~\cite{Dupacova2003}, Heitsch and Römisch \cite{Heitsch2007}, and Römisch and Wets \cite{Roemisch2007}, have laid the mathematical foundation for relating the distance between probability distributions in scenario reduction to the stability of stochastic programming problems.
More recent literature has focused on incorporating information from the optimization problem into the scenario reduction process, e.g., by first solving an expected value problem~\cite{Morales2009}, single-scenario deterministic problems~\cite{Bruninx2016}, or by additionally including cross evaluations~\cite {Bertsimas2023}. 
Here, we refer to these existing techniques as \textbf{time-invariant scenario reduction} methods, which usually focus on applications with short-time horizons, e.g., a single time step in the capacity expansion problem in \cite{Bertsimas2023}, or up to 24 hours in the transmission expansion problem in \cite{Park2019} or the unit commitment problems in \cite{Bruninx2016}, but potentially many scenarios.
For large-scale long-term planning ESOMs spanning multiple years, on the other hand, there may be only around 40 historical weather years available, and it may only be possible to include a few of them due to the immense model complexity (e.g., 16 scenarios in \cite{Goeke2024}). 
This may render existing time-invariant scenario reduction methods inefficient, as they cannot adjust their marginal distributions over time and are limited by the small number of available scenarios.

\begin{figure}
    \centering
    \definecolor{cbBlue}{RGB}{86,180,233}
\definecolor{cbGreen}{RGB}{0,158,60}
\definecolor{cbYellow}{RGB}{255,242,120}

\begin{tikzpicture}[>={Triangle[length=3pt,width=3pt]}]

\def\cellw{0.55}
\def\cellh{0.55}
\def\Nt{6}

% --------------------
% FIGURE 1
% --------------------
\begin{scope}[shift={(-1.4,0)}]
  \coordinate (fig1) at (0,0);
  \foreach \i in {1,...,\Nt} {
    \node[font=\small] at ({\i*\cellw - 0.5*\cellw},{3*\cellh + 0.22}) {$t_{\i}$};
  }
  \node[anchor=east,font=\small] at (0.0, {2.5*\cellh}) {$\xi_1$};
  \node[anchor=east,font=\small] at (0.0, {1.5*\cellh}) {$\xi_2$};
  \node[anchor=east,font=\small] at (0.0, {0.5*\cellh}) {$\xi_3$};
  \foreach \i in {1,...,\Nt} {
    \pgfmathsetmacro{\x}{\i*\cellw}
    \fill[cbBlue!50]   ({\x-\cellw},{2*\cellh}) rectangle ({\x},{3*\cellh});
    \draw[line width=0.3pt] ({\x-\cellw},{2*\cellh}) rectangle ({\x},{3*\cellh});
    \fill[cbGreen!40]  ({\x-\cellw},{\cellh})   rectangle ({\x},{2*\cellh});
    \draw[line width=0.3pt] ({\x-\cellw},{\cellh})   rectangle ({\x},{2*\cellh});
    \fill[cbYellow!100] ({\x-\cellw},{0})       rectangle ({\x},{\cellh});
    \draw[line width=0.3pt] ({\x-\cellw},{0})       rectangle ({\x},{\cellh});
  }
  \draw[line width=0.6pt] (0,0) rectangle ({\Nt*\cellw},{3*\cellh});
  \coordinate (fig1east)   at ({\Nt*\cellw},{1.5*\cellh});
  \coordinate (fig1south)  at ({0.5*\Nt*\cellw},0);
\end{scope}
\node[below=0mm of fig1south,font=\small] {Original sample};

% --------------------
% FIGURE 2
% --------------------
\begin{scope}[shift={(3.2,1.375)}]
  \foreach \i in {1,...,\Nt} {
    \node[font=\small] at ({\i*\cellw - 0.5*\cellw},{2*\cellh + 0.22}) {$t_{\i}$};
  }
  \node[anchor=west,font=\small] (fig2lbl1) at (-0.55, {1.5*\cellh}) {$\zeta^{\scriptscriptstyle\mathrm{I}}_1$};
  \node[anchor=west,font=\small] (fig2lbl2) at (-0.55, {0.5*\cellh}) {$\zeta^{\scriptscriptstyle\mathrm{I}}_2$};
  \foreach \i in {1,...,\Nt} {
    \pgfmathsetmacro{\x}{\i*\cellw}
    \fill[cbYellow!100] ({\x-\cellw},{\cellh}) rectangle ({\x},{2*\cellh});
    \begin{scope}
      \clip ({\x-\cellw},{\cellh}) rectangle ({\x},{2*\cellh});
      \foreach \j in {-2,...,3} {
        \draw[cbBlue!70, line width=1.8pt] ({\x-\cellw + \j*\cellw/3},{\cellh}) -- ++(1,1);
      }
    \end{scope}
    \fill[cbGreen!40] ({\x-\cellw},{0}) rectangle ({\x},{\cellh});
    \draw[line width=0.3pt] ({\x-\cellw},{\cellh}) rectangle ({\x},{2*\cellh});
    \draw[line width=0.3pt] ({\x-\cellw},{0})      rectangle ({\x},{\cellh});
  }
  \draw[line width=0.6pt] (0,0) rectangle ({\Nt*\cellw},{2*\cellh});
  \coordinate (fig2south)  at ({0.5*\Nt*\cellw},0);
\end{scope}
\node[below=0mm of fig2south,font=\small] {Time-invariant SR};

% --------------------
% FIGURE 3
% --------------------
\begin{scope}[shift={(3.2,-0.825)}]
  \foreach \i in {1,...,\Nt} {
    \node[font=\small] at ({\i*\cellw - 0.5*\cellw},{2*\cellh + 0.22}) {$t_{\i}$};
  }
  \node[anchor=west,font=\small] (fig3lbl1) at (-0.55, {1.5*\cellh}) {$\zeta^{\scriptscriptstyle\mathrm{V}}_1$};
  \node[anchor=west,font=\small] (fig3lbl2) at (-0.55, {0.5*\cellh}) {$\zeta^{\scriptscriptstyle\mathrm{V}}_2$};
  \foreach \i in {1,...,\Nt} {
    \pgfmathsetmacro{\x}{\i*\cellw}
    \ifnum\i<3
      \fill[cbYellow!100] ({\x-\cellw},{\cellh}) rectangle ({\x},{2*\cellh});
      \begin{scope}
        \clip ({\x-\cellw},{\cellh}) rectangle ({\x},{2*\cellh});
        \foreach \j in {-2,...,3} {
          \draw[cbBlue!70, line width=1.8pt] ({\x-\cellw + \j*\cellw/3},{\cellh}) -- ++(1,1);
        }
      \end{scope}
    \else\ifnum\i<6
      \fill[cbBlue!50] ({\x-\cellw},{\cellh}) rectangle ({\x},{2*\cellh});
    \else
      \fill[cbYellow!100] ({\x-\cellw},{\cellh}) rectangle ({\x},{2*\cellh});
    \fi\fi
    \ifnum\i<3
      \fill[cbGreen!40] ({\x-\cellw},{0}) rectangle ({\x},{\cellh});
    \else\ifnum\i<6
      \fill[cbYellow!100] ({\x-\cellw},{0}) rectangle ({\x},{\cellh});
      \begin{scope}
        \clip ({\x-\cellw},{0}) rectangle ({\x},{\cellh});
        \foreach \j in {-2,...,3} {
          \draw[cbGreen!55, line width=1.8pt] ({\x-\cellw + \j*\cellw/3},{0}) -- ++(1,1);
        }
      \end{scope}
    \else
      \fill[cbGreen!45] ({\x-\cellw},{0}) rectangle ({\x},{\cellh});
      \begin{scope}
        \clip ({\x-\cellw},{0}) rectangle ({\x},{\cellh});
        \foreach \j in {-2,...,3} {
          \draw[cbBlue!85, line width=1.8pt] ({\x-\cellw + \j*\cellw/3},{0}) -- ++(1,1);
        }
      \end{scope}
    \fi\fi
    \draw[line width=0.3pt] ({\x-\cellw},{\cellh}) rectangle ({\x},{2*\cellh});
    \draw[line width=0.3pt] ({\x-\cellw},{0})      rectangle ({\x},{\cellh});
  }
  \draw[line width=0.6pt] (0,0) rectangle ({\Nt*\cellw},{2*\cellh});
  \coordinate (fig3south)  at ({0.5*\Nt*\cellw},0);
\end{scope}
\node[below=0mm of fig3south,font=\small] {Time-variant SR};

\draw[->, thick]
  (fig1east) .. controls +(0.5,0) and +(-0.5,0) .. ([xshift=2pt]$(fig2lbl1.west)!0.5!(fig2lbl2.west)$);
\draw[->, thick]
  (fig1east) .. controls +(0.5,0) and +(-0.5,0) .. ([xshift=2pt]$(fig3lbl1.west)!0.5!(fig3lbl2.west)$);

\end{tikzpicture}
    \vspace{-0.7cm}
    \caption{Schematic representation of the conventional time-invariant (I) and the proposed time-variant (V) scenario reduction (SR) methods with three scenarios and six time steps.}
    \label{fig:SchematicSR}
    \vspace{-0.15cm}
\end{figure}
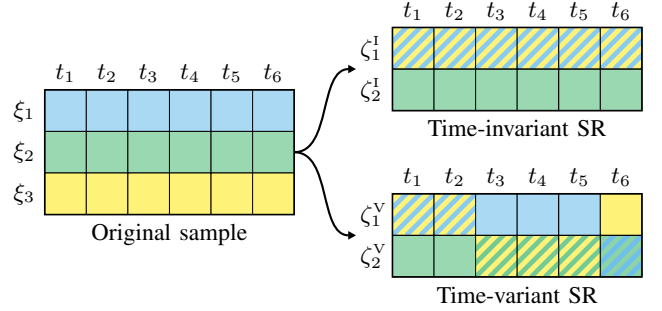

At the same time, it facilitates a new branch of research targeting the reduction of a small sample of scenarios with long time horizons by breaking up the temporal chronology of the underlying stochastic process and focusing on capturing the time-index marginal distributions.
In this letter, we propose a novel design that allows scenario aggregations to change over time, which we term \textbf{time-variant scenario reduction}.
Figure~\ref{fig:SchematicSR} illustrates the conventional time-invariant (I) and the proposed time-variant (V) scenario reduction (SR) frameworks using a schematic example with six time steps $t_1 - t_6$, aggregating three scenarios, $\xi_1 - \xi_3$, into two, $\zeta_1$ and $\zeta_2$.
These changing scenario aggregations further cause scenario probabilities to vary over time.
Let $i \in \mathcal{I}$ and $j \in \mathcal{J}$ be the original and the reduced set of scenarios, respectively, and $t \in \mathcal{T} = \{1,2,...,T\}$ the set of time steps. Furthermore, we use $\mathcal{I}_{j,t} \subseteq \mathcal{I}$ to denote the subsets of scenarios that get aggregated to scenario $j$ in time step $t \in \mathcal{T}$. Then we define the time-variant probability $\pi_{t,j}$ as:
\begin{equation}\label{eq:TV_prob_def}
    \pi_{t,j} = \sum_{i \in \mathcal{I}_{j,t}} \pi_i, \forall j \in \mathcal{J}, \forall t \in \mathcal{T}.
\end{equation}
As we show with an illustrative proof of concept, this unlocks enormous flexibility to reduce model complexity while still accurately capturing the realizations and probabilities of critical scenarios that determine first-stage investment decisions when using the original distribution.

In the remainder of this letter, Section~\ref{sec:PoC} introduces a proof of concept demonstrating the potential of the proposed time-variant scenario reduction method, Section~\ref{sec:results} presents and discusses the results obtained, and Section~\ref{sec:conclusions} concludes.
\section{Proof of concept}
\label{sec:PoC}

\subsection{Towards time-variant scenario reduction}
\label{sec:ProbForecast}
For illustrative purposes, we consider a joint probabilistic forecast with 20 equiprobable scenarios for photovoltaic (PV) and wind power production over 96 hours (four days)~\cite{Lueth2024}, which we must reduce to solve the two-stage stochastic ESOM introduced in Section~\ref{sec:GEP}. Here, we utilize the popular k-means clustering algorithm \cite{Lloyd1982}.

Let vector $w_{t,i}$ denote the stacked PV and wind power production forecast in time step~$t$ and scenario~$i$. For the conventional time-invariant SR, k-means can be applied straightforwardly to the set of scenarios spanning the joint realizations over the full time horizon $w_i^\intercal =  [w_{1,i}^\intercal,\dots,w_{|\mathcal{T}|,i}^\intercal]$.
For the novel time-variant SR, we first break the temporal chronology of the scenarios and apply k-means clustering to $w_{t,i}$ for each time step $t \in \mathcal{T}$ individually. Afterward, to restore temporally chronological time series for scenarios $j$ in the reduced set~$\mathcal{J}$, we solve, for each time step $t \in \{1,\dots, T-1\}$, the following linear intertemporal assignment problem with weight $x_{j,j',t} \in \{0,1\}, \forall j,j' \in \mathcal{J}$ and cost $c_{j,j',t} = \lVert w_{t+1,j} - w_{t,j'} \rVert_2^2$:
\begin{align}
    \min_{x_{j,j',t}} ~ & \sum_{j \in \mathcal{J}} \sum_{j' \in \mathcal{J}} x_{j,j',t} \cdot c_{j,j',t}, \\
    \mathrm{s.t.} ~ & \sum_{j' \in \mathcal{J}} x_{j,j',t} = 1, \forall j \in \mathcal{J}, \sum_{j \in \mathcal{J}} x_{j,j',t} = 1, \forall j' \in \mathcal{J}.
\end{align}
This procedure yields subsets $\mathcal{I}_{j,t} \subseteq \mathcal{I}$ of scenarios in $\mathcal{I}$ that get aggregated to scenario $j$ in time step $t$, which we use to derive the time-variant probabilities $\pi_{t,j}$ according to Equation~\ref{eq:TV_prob_def}.
For illustration, we show the probabilities for both methods applied to the probabilistic forecast at hand across a set of selected time steps and four scenarios in the reduced set in Figure~\ref{fig:probabilities}.
In stark contrast to the time-invariant SR (a) that constantly allocates large probability mass to two scenarios, the time-variant SR (b) flexibly adjusts scenario probabilities over time, allowing it to accurately capture scenario probabilities at critical time steps that drive first-stage investment decisions.

\subsection{Stochastic generation expansion planning problem}
\label{sec:GEP}

We now introduce the scenario-based two-stage stochastic generation expansion planning problem $\mathbf{SGEP}$ in extensive form, using lowercase and uppercase letters to denote decision variables and parameters, respectively.
Random variables are highlighted using a hat, and optimal values of decision variables are marked with an asterisk.

\begin{figure}
    \centering
    \begin{tikzpicture}[scale=1]
% -------------------------------------------------------
% Parameters
% -------------------------------------------------------
\def\cellw{0.5}
\def\cellh{0.5}
\def\firstTime{15}
\def\lastTime{25}
\def\gap{0.4}        % gap between TI column and TV block
% -------------------------------------------------------
% Read CSV file (time_step, cluster_1..cluster_4) for TV case
% -------------------------------------------------------
\pgfplotstableread[col sep=comma]{tikz/TV_apriori_global_cluster_probabilities.csv}\probtable
% -------------------------------------------------------
% Time-invariant probabilities (clusters 1..4)
% -------------------------------------------------------
\def\tiA{0.5}
\def\tiB{0.05}
\def\tiC{0.15}
\def\tiD{0.3}
% -------------------------------------------------------
% Calculate size
% -------------------------------------------------------
\pgfmathtruncatemacro{\numcols}{\lastTime-\firstTime+1}
\def\numrows{4}
% -------------------------------------------------------
% Color setup: viridis
% -------------------------------------------------------
\pgfplotsset{colormap/viridis}
\def\colormin{0}
\def\colormax{0.6}
% -------------------------------------------------------
% TI column row labels (left of TI column)
% -------------------------------------------------------
\foreach \c [count=\i from 1] in {1,2,3,4} {
  \node[anchor=east] at (-0.01, {(\numrows-\i+0.5)*\cellh}) {\small $\zeta^{\mathrm{I}}_{\c}$};
}
% -------------------------------------------------------
% Draw TI column (single column, x in [0, \cellw])
% -------------------------------------------------------
\foreach \p [count=\i from 1] in {\tiA,\tiB,\tiC,\tiD} {
  \pgfmathsetmacro{\metaval}{(\p-\colormin)/(\colormax-\colormin)*1000}
  \pgfplotscolormapdefinemappedcolor{\metaval}
  \definecolor{cellcolor}{rgb}{\pgfmathresult}
  \pgfmathsetmacro{\y}{(\numrows-\i)*\cellh}
  \fill[cellcolor] (0,\y) rectangle ++(\cellw,\cellh);
  \draw (0,\y) rectangle ++(\cellw,\cellh);
}
\draw[very thick] (0,0) rectangle (\cellw,\numrows*\cellh);
% Panel label (a) under TI
\node[anchor=north,font=\small] at ({0.5*\cellw}, -0.15) {(a)};
% -------------------------------------------------------
% TV block: shift right by \cellw + \gap
% -------------------------------------------------------
\begin{scope}[xshift={(\cellw+\gap)*1.75cm}]
  % Column headers
  \foreach \t [count=\j from 1] in {\firstTime,...,\lastTime} {
    \node at (\j*\cellw - 0.5*\cellw, {\numrows*\cellh + 0.2}) {\small $t_{\t}$};
  }
  % Row labels (clusters) for TV
  \foreach \c [count=\i from 1] in {1,2,3,4} {
    \node[anchor=east] at (-0.01, {(\numrows-\i+0.5)*\cellh}) {\small $\zeta^{\mathrm{V}}_{\c}$};
  }
  % Draw cells
  \foreach \t [count=\j from 1] in {\firstTime,...,\lastTime} {
    \pgfplotstablegetelem{\t-1}{cluster_1}\of{\probtable}\edef\pA{\pgfplotsretval}
    \pgfplotstablegetelem{\t-1}{cluster_2}\of{\probtable}\edef\pB{\pgfplotsretval}
    \pgfplotstablegetelem{\t-1}{cluster_3}\of{\probtable}\edef\pC{\pgfplotsretval}
    \pgfplotstablegetelem{\t-1}{cluster_4}\of{\probtable}\edef\pD{\pgfplotsretval}
    \foreach \p [count=\i from 1] in {\pA,\pB,\pC,\pD} {
      \pgfmathsetmacro{\metaval}{(\p-\colormin)/(\colormax-\colormin)*1000}
      \pgfplotscolormapdefinemappedcolor{\metaval}
      \definecolor{cellcolor}{rgb}{\pgfmathresult}
      \pgfmathsetmacro{\x}{(\j-1)*\cellw}
      \pgfmathsetmacro{\y}{(\numrows-\i)*\cellh}
      \fill[cellcolor] (\x,\y) rectangle ++(\cellw,\cellh);
      \draw (\x,\y) rectangle ++(\cellw,\cellh);
    }
  }
  % Outer border
  \draw[very thick] (0,0) rectangle (\numcols*\cellw,\numrows*\cellh);
  % Panel label (b) under TV
  \node[anchor=north,font=\small] at ({0.5*\numcols*\cellw}, -0.15) {(b)};
\end{scope}
% -------------------------------------------------------
% Colorbar on the right (shifted further to account for TI column + gap)
% -------------------------------------------------------
\begin{scope}[xshift={(\cellw+\gap+\numcols*\cellw)+7.2cm}]
  \foreach \k in {0,...,10} {
    \pgfmathsetmacro{\metaval}{\k*100}
    \pgfplotscolormapdefinemappedcolor{\metaval}
    \definecolor{barcolor}{rgb}{\pgfmathresult}
    \pgfmathsetmacro{\ypos}{\k*0.18}
    \fill[barcolor] (0, \ypos) rectangle (0.3, \ypos+0.18);
  }
  \draw (0,0) rectangle (0.3,1.98);
  \node[anchor=west,font=\small] at (0.275,0) {0};
  \node[anchor=west,font=\small] at (0.275,1.98) {0.6};
\end{scope}
\end{tikzpicture}
    \vspace{-0.75cm}
    \caption{(a) Time-invariant and (b) time-variant probabilities in the reduced scenario sets for selected time steps.}
    \label{fig:probabilities}
    \vspace{-0.15cm}
\end{figure}
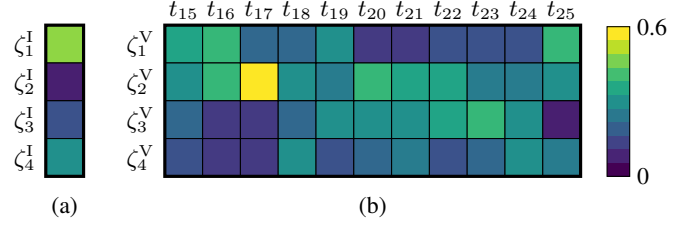

In the first stage, i.e., before the uncertainty is realized, investment in battery storage capacity $s$ and thermal generation capacity $x$ at costs $Z^{\mathrm{inv}}$ and $C^{\mathrm{inv}}$, respectively, must be decided. In the second stage, i.e., after uncertainty is realized, for every time step $t \in \mathcal{T}$, the inelastic electricity demand $D_t$ must be met for any realization $\omega \in \Omega$ of uncertain power production from PV, $\widehat{J}_{t,\omega}$, and wind, $\widehat{F}_{t,\omega}$, utilizing the existing thermal power plants, $q_{t,\omega}$, at cost $C^{\mathrm{th}}$, and the newly invested thermal power plants $g_{t,\omega}$ at cost $C^{\mathrm{op}}$. Electrical load shedding $d^{\mathrm{sh}}_{t,\omega}$ incurs high penalty costs $C^{\mathrm{sh}}$.
The objective of the $\mathbf{SGEP}$ problem is to minimize the sum of first-stage investment and expected second-stage operational costs:
\begin{equation}\label{eq:obj}
\begin{split}
    \min \quad {} & =
        C^{\mathrm{inv}} x +  Z^{\mathrm{inv}} s
        \vphantom{\sum_{\omega \in \Omega}}
    + \\ %[2pt]
    &
        \sum_{\omega \in \Omega}
        \left[ 
            \sum_{t \in \mathcal{T}}
             \pi_{t,\omega}
            \left(
                C^{\mathrm{th}} q_{t,\omega} 
                + C^{\mathrm{op}} g_{t,\omega} 
                + C^{\mathrm{sh}} d^{\mathrm{sh}}_{t,\omega}
            \right)
        \right].
\end{split}
\end{equation}
where $\pi_{t,\omega} \in [0,1]$ is the probability of scenario $\omega$ in time step~$t$, shown in Figure~\ref{fig:probabilities}. Note that for the time-invariant SR, this probability is constant, such that $\pi_{t,\omega} = \pi_{\omega}, \forall t \in \mathcal{T}$.

First-stage investment decisions $x$ and $s$ are bounded above by sufficiently large limits $X^{\mathrm{max}}$ and $S^{\mathrm{max}}$, respectively:

\vspace{5pt}
% \vspace{-3pt}
\begin{subequations}\label{eq:FirstStageConstraints}
\noindent
\begin{minipage}{.5\linewidth}
\begin{equation}
0 \leq x \leq X^{\mathrm{max}}, \label{eq:InvThLimit}
\end{equation}
\end{minipage}%
\begin{minipage}{.5\linewidth}
\begin{equation}
0 \leq s \leq S^{\mathrm{max}}. \label{eq:InvStorLimit}
\end{equation}
\end{minipage}
\end{subequations}
\vspace{1pt}
\vspace{-3pt}

The second-stage operational constraints for every time step $t \in \mathcal{T}$ and scenario $\omega \in \Omega$, are given by:

\vspace{5pt}
\begin{subequations}\label{eq:SecondStageConstraints}
\noindent
\begin{minipage}{.5\linewidth}
\hfill\makebox[0.5\linewidth][l]{$\displaystyle 0 \leq g_{t,\omega} \leq x,$}\hfill\refstepcounter{equation}(\theequation)\label{eq:InvTh}
\end{minipage}%
\begin{minipage}{.5\linewidth}
\hfill\makebox[0.5\linewidth][l]{$\displaystyle 0 \leq q_{t,\omega} \leq Q^{\mathrm{max}},$}\hfill\refstepcounter{equation}(\theequation)\label{eq:ExistProdLimit}
\end{minipage}\\[6pt]
\noindent
\begin{minipage}{.5\linewidth}
\hfill\makebox[0.5\linewidth][l]{$\displaystyle 0 \leq f^{\mathrm{sp}}_{t,\omega} \leq \widehat{F}_{t,\omega},$}\hfill\refstepcounter{equation}(\theequation)\label{eq:SpillWind}
\end{minipage}%
\begin{minipage}{.5\linewidth}
\hfill\makebox[0.5\linewidth][l]{$\displaystyle 0 \leq j^{\mathrm{sp}}_{t,\omega} \leq \widehat{J}_{t,\omega},$}\hfill\refstepcounter{equation}(\theequation)\label{eq:SpillPV}
\end{minipage}\\[6pt]
\noindent
\begin{minipage}{.5\linewidth}
\hfill\makebox[0.5\linewidth][l]{$\displaystyle 0 \leq p^{\mathrm{ch}}_{t,\omega}, p^{\mathrm{dis}}_{t,\omega} \leq s,$}\hfill\refstepcounter{equation}(\theequation)\label{eq:InvStor}
\end{minipage}%
\begin{minipage}{.5\linewidth}
\hfill\makebox[0.5\linewidth][l]{$\displaystyle 0 \leq e_{t,\omega} \leq 2 \cdot s,$}\hfill\refstepcounter{equation}(\theequation)\label{eq:InvStorEnergy}
\end{minipage}
\vspace{-14pt}
\begin{align}
& e_{t,\omega} = (1-\gamma) \cdot e_{t-1,\omega} + \eta \cdot p^{\mathrm{ch}}_{t,\omega} - 1/\eta \cdot p^{\mathrm{dis}}_{t,\omega}, \label{eq:StorInventory} \\
\begin{split}
& q_{t,\omega} + g_{t,\omega} + p^{\mathrm{dis}}_{t,\omega} - p^{\mathrm{ch}}_{t,\omega} = \\
& \quad D_t - d^{\mathrm{sh}}_{t,\omega} - \widehat{F}_{t,\omega} + f^{\mathrm{sp}}_{t,\omega} - \widehat{J}_{t,\omega} + j^{\mathrm{sp}}_{t,\omega}.
\end{split} \label{eq:ED_power_bal}
\end{align}
\end{subequations}
% \vspace{2pt}

Constraints~\eqref{eq:InvTh} and \eqref{eq:ExistProdLimit} limit the power production of new and existing thermal units to $x$ and $Q^{\mathrm{max}}$, respectively. Uncertain wind and PV power production can be fully spilled at no cost by $f^{\mathrm{sp}}_{t,\omega}$ and $j^{\mathrm{sp}}_{t,\omega}$, respectively, given Constraints~\eqref{eq:SpillWind} and \eqref{eq:SpillPV}. Storage charging $p^{\mathrm{ch}}_{t,\omega}$ and discharging $p^{\mathrm{dis}}_{t,\omega}$ is limited by Constraints~\eqref{eq:InvStor}, whereas Constraints~\eqref{eq:InvStorEnergy} bound the stored energy $e_{t,\omega}$ between zero and the invested capacity $s$, times two (by assumption). The energy stored is determined by inventory constraint~\eqref{eq:StorInventory} considering an intertemporal loss rate~$\gamma$ and a symmetric dis-/charge efficiency~$\eta$ and assuming that the initial storage level is equal to zero in all scenarios, i.e., $e_{0,\omega} = 0, \forall \omega \in \Omega$. Finally, the power balance in each scenario and time step is ensured by Constraints~\eqref{eq:ED_power_bal}.

\section{Numerical results}
\label{sec:results}

\begin{figure}
    \centering
    % \documentclass[border=3pt]{standalone}
% \usepackage{pgfplots}
% \pgfplotsset{compat=1.18}
% \pgfplotsset{table/col sep=comma}   % data file is comma-separated
% \usepackage{xcolor}
\definecolor{colorI}{RGB}{240,150,30}   % time-invariant (orange)
\definecolor{colorV}{RGB}{30,144,255}   % time-variant (blue)

% \begin{document}
 \begin{tikzpicture}
  \begin{axis}[
      width=\columnwidth, height=5.5cm,
      xlabel={Size of reduced scenario set $\mathcal{J}$}, ylabel={RAE [\%]},
      xmin=1, xmax=10, ymin=0, ymax=100,
      xtick distance=1,
      axis lines=left, axis line style={-},
      tick align=outside, grid=none,
      ylabel style={yshift=-4pt},
      label style={font=\small},
      tick label style={font=\small},
      legend style={font=\small, draw=black, fill=white},
      legend pos=north east, legend cell align=left,
  ]
  \foreach \m in {Jan,Feb,Mar,Apr,May,Jun,Jul,Aug,Sep,Oct,Nov,Dec}{%
      \addplot[colorI, opacity=0.22, line width=1.1pt, forget plot]
          table[col sep=comma, x=clusters, y=apr_\m]{Figures/spaghetti/spaghetti_base_data.csv};}
  \foreach \m in {Jan,Feb,Mar,Apr,May,Jun,Jul,Aug,Sep,Oct,Nov,Dec}{%
      \addplot[colorV, opacity=0.22, line width=1.1pt, forget plot]
          table[col sep=comma, x=clusters, y=tv_\m]{Figures/spaghetti/spaghetti_base_data.csv};}

      % --- bold means (thick); shown in the plot, explained via the legend note ---
      \addplot[colorI, line width=1.3pt, forget plot] table[col sep=comma, x=clusters, y=apr_mean]{Figures/spaghetti/spaghetti_base_data.csv};
      \addplot[colorV, line width=1.3pt, forget plot] table[col sep=comma, x=clusters, y=tv_mean]{Figures/spaghetti/spaghetti_base_data.csv};
      % --- legend: one faint method line each (matches the individual windows) + mean note ---
      \addlegendimage{colorI, opacity=0.5, line width=1.1pt}
      \addlegendentry{Time-invariant}
      \addlegendimage{colorI, line width=1.3pt}
      \addlegendentry{Time-invariant (mean)}
      \addlegendimage{colorV, opacity=0.5, line width=1.1pt}
      \addlegendentry{Time-variant}
      \addlegendimage{colorV, line width=1.3pt}
      \addlegendentry{Time-variant (mean)}
      
  \end{axis}
  \end{tikzpicture}
% \end{document}
    \vspace{-0.35cm}
    \caption{Relative approximation error (RAE) for different time windows and numbers of scenarios in the reduced set.}
    \label{fig:RAE_sensitivity}
\end{figure}
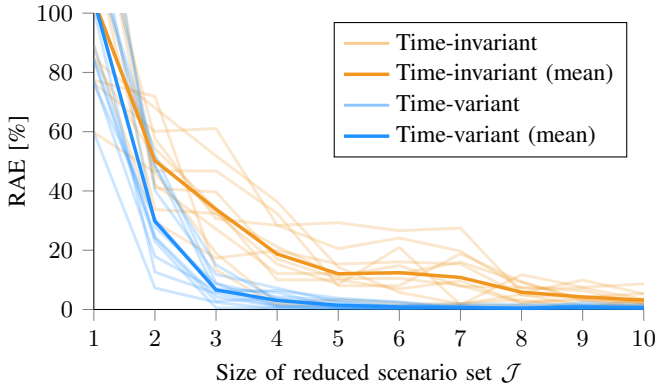

We solve the $\mathbf{SGEP}$ problem~\eqref{eq:obj}--\eqref{eq:SecondStageConstraints} using the time-invariant and time-variant SR presented in Section~\ref{sec:ProbForecast}. To evaluate the solution quality, we fix the first-stage investment decisions obtained and run independent operational problems for all 20 scenarios in the original sample. Afterward, we calculate the resulting expected investment and operational costs $W$, and derive the relative approximation error $\mathrm{RAE} = (W^{r}-W^{\mathrm{S}})/W^{\mathrm{S}} \cdot 100\%$ (cf. \cite{Pflug2001}), where $r$ is either time-invariant ($\mathrm{I}$) or time-variant ($\mathrm{V}$), and $W^{\mathrm{S}}$ the optimal solution to $\mathbf{SGEP}$ for the original sample.
Figure~\ref{fig:RAE_sensitivity} shows the $\mathrm{RAE}$ for both SR frameworks, varying the number of scenarios in the reduced set and the 4-day time window of the original sample (one per month).
The plot indicates that the time-variant SR converges to the solution for the full scenario set much faster and more robustly, while achieving good solution accuracy from as few as four scenarios in the reduced set.

Table~\ref{tab:investment_comparison_avg} shows detailed optimization model results, including the expected energy not served ($\mathrm{EENS}$) for $|\mathcal{J}| = 4$ scenarios, averaged over all samples.
Compared with the results of the full original sample with 20 scenarios, the time-invariant (I) SR method with four scenarios performs poorly: it substantially underestimates the optimal investment in thermal generation capacity, $x^*$, while overestimating the investment in storage, $s^*$. This leads to a massive increase in $\mathrm{EENS}$ and an RAE of around $\qty{19}{\percent}$.
In contrast, the time-variant SR method accurately reproduces the optimal investments in both thermal generation and storage with only four scenarios, achieving an $\mathrm{RAE}$ of just $\qty{3}{\percent}$. This reduction in the number of scenarios corresponds to an approximate $\qty{80}{\percent}$ decrease in model size, as measured by the number of variables. To get the same performance with the time-invariant SR used here, you would need at least 11 scenarios in the reduced set.

The superior performance of the time-variant SR stems from its ability to accurately capture scenario realizations and probabilities at the time steps that drive investment decisions by adjusting time-specific marginal distributions. This ability is fundamentally unavailable in traditional time-invariant SR methods because scenario probabilities must remain constant over time, as shown in Figure~\ref{fig:probabilities}.
  
\begin{table}%[htbp]
\centering
\caption{$\mathbf{SGEP}$ results for the full sample, time-invariant~(I) and time-variant (V) SR methods with $|\mathcal{J}| = 4$ scenarios averaged over the different time windows.}
\label{tab:investment_comparison_avg}
\small
\begin{tabular}{lccccccc}
\toprule
&  $x^*$ & $s^*$ & $W$ & $\mathrm{RAE}$ & $\mathrm{EENS}$ \\
&  (\unit{\mega\watt}) & (\unit{\mega\watt}) & (\unit{\kilo\EUR}) & (\%) & (\unit{\mega\watt\hour}) \\
\midrule
Full & 5054.4 & 2549.3 & 102.4 & 0.00 & 113.9 \\
I SR & 4217.2 & 2797.3 & 121.5 & 18.71 & 2815.8 \\
V SR & 4497.7 & 3279.1 & 105.7 & 3.06 & 712.9 \\
\bottomrule
\end{tabular}
\end{table}
  
\section{Summary and conclusions}
\label{sec:conclusions}
In this letter, we propose a novel time-variant scenario reduction framework that leverages the temporal nature of scenarios in long-term planning problems by allowing scenario aggregations to vary over time. This introduces time-variant scenario probabilities, enabling accurate capture of scenario realizations and probabilities at the time steps that drive investment decisions. We demonstrate the capabilities of the proposed approach using an illustrative two-stage stochastic generation expansion planning problem with uncertain wind and PV production. Future research should focus on deriving appropriate time-variant scenario distance metrics and algorithms to avoid the two-step procedure used in this paper. Furthermore, it should be investigated how mathematical stability results can be extended to the time-variant case to derive performance guarantees.

\section*{Acknowledgement}
The work of J. M. Morales and S. Pineda is supported by the Spanish Ministry of Science, Innovation and Universities through project PID2023-148291NB-I00 and by the Department of Universities, Research and Innovation of the Regional Government of Andalusia through project DGP\_PIDI\_2024\_00851.

\bibliographystyle{IEEEtran}
\bibliography{paper/TVSR}

\end{document}